\documentclass[leqno,12pt]{article}

\usepackage{latexsym}
 \usepackage{graphicx}
\usepackage{amsmath}
\usepackage{amssymb}
\usepackage{cite}
\usepackage{color}

\usepackage[margin=1.3in, top=1.3in, bottom=1.3in]{geometry}

\newcommand{\nc}{\newcommand}
\nc{\nt}{\newtheorem}
\nt{thm}{Theorem}[section]
\nt{cor}[thm]{Corollary}
\nt{prop}[thm]{Proposition}
\nt{lem}[thm]{Lemma}
\nt{defn}[thm]{Definition}
\nt{rem}[thm]{Remark}
\nt{exa}[thm]{Example}
\nt{ass}[thm]{Assumption}
\nc{\ip}[2]{\mbox{$\langle #1,#2 \rangle$}}
\nc{\pf}{\noindent{\bf Proof\ \ }}
\nc{\finpf}{\hfill{$\Box$}\linespace}
\nc{\linespace}{\vspace{\baselineskip} \noindent}
\nc{\R}{{\bf R}}
\nc{\bR}{{\overline{\R}}}
\nc{\C}{{\bf C}}
\nc{\E}{{\bf E}}
\nc{\Y}{{\bf Y}}
\nc{\RE}{\mbox{\rm Re}\,}
\nc{\Rn}{{\bf R}^n}
\nc{\Mn}{{\bf M}^n}
\nc{\bx}{\bar{x}}
\nc{\by}{\bar{y}}
\nc{\e}{\epsilon}
\nc{\inT}{\mbox{\rm int}\,}
\nc{\cl}{\mbox{\rm cl}\,}
\nc{\gph}{\mbox{\rm gph}\,}
\nc{\conv}{\mbox{\rm conv}\,}
\nc{\rt}{\rightarrow}
\nc{\tra}{\mbox{\rm tr}\,}
\def\ck{{\mathcal K}}

\def\slof{\overline{|\nabla f|}}

\nc{\xbar}{{\overline{x}}}
\nc{\vbar}{{\overline{v}}}
\nc{\yhat}{{\widehat{y}}}
\nc{\xhat}{{\widehat{x}}}

\nc{\und}{\quad\mbox{ and }\quad}
\nc{\emp}{\ensuremath{\varnothing}}

\def\tto{\;{\lower 1pt \hbox{$\rightarrow$}}\kern -12pt
           \hbox{\raise 2.8pt \hbox{$\rightarrow$}}\;}
\newenvironment{myequation}{\setcounter{equation}{\value{thm}}
   \begin{equation}}{\addtocounter{thm}{1}\end{equation}}

\nc{\bmye}{\begin{myequation}}
\nc{\emye}{\end{myequation}}

\begin{document}
\title{
Transversality and alternating projections for nonconvex sets
}
\author{
D. Drusvyatskiy\thanks{Department of Mathematics, University of Washington, Seattle, WA 98195;~~~~~~~\mbox{}  \texttt{www.math.washington.edu/$\sim$ddrusv}.
Research supported in part by AFOSR YIP FA9550-15-1-0237.}
\and
A.D. Ioffe\thanks{Mathematics Department, Technion-Israel Institute of Technology, Haifa, Israel 32000;  \texttt{ioffe@math.technion.ac.il}. Research supported in part by the US-Israel Binational Science Foundation Grant 2008261.}
\and
 A.S. Lewis\thanks{ORIE, Cornell University, Ithaca, NY 14853, U.S.A.
\texttt{people.orie.cornell.edu/aslewis}.
Research supported in part by National Science Foundation Grant DMS-1208338 and by the US-Israel Binational Science Foundation Grant 2008261.}
}
\maketitle

\begin{abstract}
We consider the method of alternating projections for finding a point in the intersection of two closed sets, possibly nonconvex.  Assuming only the standard transversality condition (or a weaker version thereof), we prove local linear convergence.  When the two sets are semi-algebraic and bounded, but not necessarily transversal, we nonetheless prove subsequence convergence.
\end{abstract}
\medskip

\medskip

\medskip
\begin{center}
Communicated by Michael Todd.
\end{center}

\section{Tranversality}
The classical idea of {\em transversality} at an intersection point of two smooth manifolds in a Euclidean space $\E$ means that the two tangent spaces at the intersection point sum to $\E$:  equivalently, their two orthogonal complements (the {\em normal} spaces) intersect trivially.  Transversality assumptions are ubiquitous throughout analysis, guaranteeing in particular some stability for intersections under small perturbations.  For two convex sets, for example, this stability requires that the intersection point lies on no separating hyperplane, a property also easy to express using normal vectors.

For more general sets, that may be neither convex nor smooth, the correct generalization of transversality involves normals rather than tangents.  Finite-dimensional variational analysis allows a unified treatment through a single, fundamental geometric construction:  the normal cone.  (The monographs \cite{Borwein-Zhu,CLSW,Mord_1,VA} are good references.)  In the space $\E$, we denote the inner product $\ip{\cdot}{\cdot}$ and the norm $
|\cdot|$.  The set of points in a set $X \subset \E$ at minimum distance from a point $y \in \E$ is the {\em projection} $P_X(y)$, and that minimum distance is $d(y,X)$.  For any point $x \in X$, vectors in the cone
\[
N^p_X(x) ~=~ \big\{ \lambda u : \lambda \in \R_+,~ u \in \E,~ x \in P_X(x+u) \big\}
\]
are called {\em proximal normals} to $X$ at $x$.  Limits of proximal normals to $X$ at points $x_n \in X$ approaching $x$ are called {\em normals}, and comprise the {\em normal cone} $N_X(x)$.  When $X$ is a smooth manifold, $N_X(x)$ is the classical normal space;  when $X$ is closed and convex, on the other hand, $N_X(x)$ consists of those linear functionals maximized over $X$ by $x$.

By analogy with the classical case of manifolds, we say that two closed sets $X,Y \subset \E$ intersect {\em transversally} at a point $z \in X \cap Y$ when
\bmye \label{transversality}
N_X(z) \cap -N_Y(z) = \{0\}.
\emye
Central to many developments of variational analysis, this property is the pre-eminent theme, in particular, in \cite{Mord_1,Mord_2}.  The terminology, highlighting a fresh understanding of the classical analogy, is more recent:  an early example is \cite[p.~99]{CLSW}.

In our entirely finite-dimensional setting, this definition of transversality is equivalent to several {\em metric regularity} properties \cite{VA}, with profound consequences for stability, sensitivity, and robustness (and for algorithms \cite{L-ICM}).  Metric regularity of set-valued mappings is not our focus, but standard ``coderivative'' computations \cite{ioffe-survey,alt_genproj} prove the equivalence of property
(\ref{transversality}) with each of the following properties.
\begin{itemize}
\item
The mapping $\Phi \colon \E^2 \tto \E$ defined by $\Phi(x,y) = x-y$ for $x \in X$ and $y \in Y$ (and empty otherwise) is metrically regular at $(z,z)$ for $0$.
\item
The mapping  $\Psi \colon \E \tto \E^2 $ defined by $\Psi(w) = (X-w) \times (Y-w)$  is metrically regular at $z$ for $(0,0)$. 
\item
There exists a constant $\tau > 0$ such that
\bmye \label{regularity}
d(w,X' \cap Y') ~\le~ \tau \big(d(w,X') + d(w,Y')\big) ~~~\mbox{for all $w \in \E$ near $z$}
\emye
and all small translations $X',Y'$ of $X,Y$.  
\end{itemize}
Transversality implies in particular the existence of a constant $\tau > 0$ such that (\ref{regularity}) holds with $X'=X$ and $Y'=Y$.  We call this weaker property {\em subtransversality}, because it corresponds to metric {\em sub\/}regularity \cite{VA} for the mappings $\Phi$ and $\Psi$.  For an extended discussion of metric regularity and transversality, see \cite{ioffe-survey}.

Transversality is not just geometrically fundamental;  it is also broad in applicability.  As we discuss below, the property is strikingly unrestrictive in the class of sets to which it potentially applies, and, in concrete settings, it holds generically.

First, note that the transversality condition (\ref{transversality}) implies little about the local structure of either of the two sets $X$ and $Y$ {\em individually}.  In particular, it may hold without either set, considered alone, being ``regular'' in the standard  single-set senses in the literature:  for example, local convexity or smoothness, prox- or Clarke regularity \cite{Borwein-Zhu,CLSW,Mord_1,VA}, or super-regularity \cite{alt_genproj}\footnote{\cite[Remark 10]{noll3} states that ``intrinsic transversality [the property, weaker than transversality, that we actually use to guarantee linear convergence] amalgamates transversality and regularity aspects''.  However, that statement concerns ``regularity'' aspects that involve {\em both} sets, like the idea of ``linear regularity'' in \cite{alt_class}}.  The epigraph 
\[
\{ (x,\tau) \in \E \times \R : \tau \ge f(x) \}
\]
 of any locally Lipschitz function $f$ intersects itself transversally at the point $(0,f(0))$, so in particular transversality is potentially applicable to any ``epi-Lipschitz'' set.  In Section \ref{intrinsic}, we extend the reach of transversality further, to sets with empty interior.

Secondly, consider the broad but concrete class of {\em semi-algebraic} sets in the space $\Rn$.  Such sets are finite unions of {\em basic} sets, each of which is a finite intersection of polynomial level sets $\{x : p(x) \le 0 \}$ and their complements \cite{BCR,Coste-semi}. Given any two closed semi-algebraic sets $X,Y \subset \R^n$, as we later prove, for almost all vectors $e \in \E$, transversality holds at  every point in the intersection of the sets $X$ and $Y-e$.

To summarize, transversality is a concise, fundamental, and widely-studied geometric property, elegantly unifying classical notions, applicable to broad classes of sets, and commonly holding in concrete settings.  Among its many powerful consequences, our aim here is to describe the first proof that transversality alone also guarantees local linear convergence for the method of alternating projections.

\section{Alternating projections}
The {\em method of alternating projections} for finding a point in the intersection of two nonempty closed sets $X,Y \subset \E$ iterates the following pair of steps:  given a current point 
$x_n \in X$,
\begin{eqnarray*}
\mbox{choose} &y_n \in	P_Y(x_n) \\
\mbox{choose} &x_{n+1} \in	P_X(y_n).
\end{eqnarray*}
This conceptually simple and widely used method for solving feasibility problems has been discovered and rediscovered by a number of authors, notably von Neumann \cite{N_proj}.  We are particularly interested in {\em linear convergence} of the sequence of iterates, by which we mean that the sequence of distance from the iterates to the limit point is bounded above by a geometric sequence.

Convergence of the method for two intersecting closed {\em convex} sets is always guaranteed \cite{conv_alt}.  Linear convergence holds whenever the relative interiors of $X$ and $Y$ intersect  \cite{GPR}, and more generally assuming just the global version of subtransversality \cite{alt_class}.

Even for nonconvex sets, the method is popular: the first linear convergence guarantees in this setting appeared in \cite{alt_man,alt_genproj}, with extensions in \cite{btheory, bapp, alt_glob} and another discussion in \cite{FW}.  In particular, the main contribution in \cite{alt_genproj} was a proof that transversality at the point $z$ suffices to guarantee linear convergence, starting sufficiently close to $z$, providing at least one of the sets $X$ and $Y$ is ``super-regular'' at $z$.  Super-regularity is a nontrivial assumption, requiring the set to be not too nonconvex, and in particular implying Clarke regularity: it opens the door to a simple geometric proof.  Super-regularity is, however, superfluous.  Using a very different proof technique --- the main topic of this paper --- the authors derived the following striking theorem.  We first discuss its provenance.

\begin{thm} \label{basic}
If two closed sets in a Euclidean space intersect transversally at a point, then the method of alternating projections, started nearby, converges linearly to a point in the intersection.
\end{thm}

The August 2013 Ph.D.\ thesis of Drusvyatskiy seems to be the first reference for this result and its proof \cite[Theorem 3.2.3]{dima-thesis}, as well as the definition of ``intrinsic'' transversality on which it relies \cite[Definition 3.2.1]{dima-thesis}.  The result (without proof) was published in August 2014, in Lewis's lecture for the International Congress of Mathematicians \cite[Theorem 2.1]{L-ICM}.  The initial version of the current work \cite{coupling} was made public and submitted in January 2014.

A new and general framework for analyzing the convergence of nonconvex alternating projections (using properties of short sequences of iterates) was submitted for publication and disseminated in December 2013 --- see \cite{noll}.  That manuscript did not show that transversality implies linear convergence (Theorem \ref{basic} above).  Although result statements in a later revision \cite[Corollary 7]{noll2} (and the final published version \cite{noll3}) do subsume Theorem \ref{basic}, they do so by explicitly relying on our core argument in \cite{coupling}\footnote{The proof of \cite[Proposition 8]{noll2} states:  ``Following entirely the argument in \cite[p.6]{coupling}\ldots''.}.   It is that argument that we present here.

In general, for two intersecting sets, no matter how close an initial point is to the intersection, the iterates generated by the method of alternating projections may not tend to the intersection. However, such pathologies rarely appear in practice.  To end, we show how our basic technique proves that, if the two sets $X$ and $Y$ are semi-algebraic and one of the sets is compact, then the method of alternating projections initiated sufficiently close to the intersection always generates iterates whose distance to the intersection tends to zero.  Further regularity guarantees convergence of the sequence of iterates \cite{noll3}.  An analogous but independent result for averaged projections appears in \cite{ABRS}.   

\section{Intrinsic transversality} \label{intrinsic}
Theorem \ref{basic} is not the last word on the local linear convergence of the method of alternating projections.  It does not cover classically well-understood cases like two distinct intersecting lines in $\R^3$, a gap discussed in some recent work \cite{FW, btheory, bapp, alt_glob, noll3}.  

The approach in \cite{alt_glob} depends, among other assumptions, on subtransversality (called there ``local linear regularity'').  Using a very different proof technique, we rely only on the following property, which we show to be intermediate in strength between transversality and subtransversality. 

\begin{defn} \label{int}
{\rm
Two closed sets $X,Y \subset \E$ are {\em intrinsically transversal} at a common point if there exists an angle  $\alpha > 0$ such that, nearby, any two points $x \in X \cap Y^c$ and 
$y \in Y \cap X^c$ cannot have difference $x-y$ simultaneously making an angle strictly less than $\alpha$ with both the cones $N_Y(y)$ and $-N_X(x)$.
}
\end{defn}

\noindent
(The normal cones may equivalently be replaced by proximal normal cones \cite{coupling}.)

While intrinsic transversality is widely applicable, we make no attempt at a detailed comparison with other variants of transversality.  Some relationships may be found in \cite{coupling}, and extensive discussions appear in \cite{noll3} and \cite{kruger,kruger-thao}.  Our focus here is rather a clear and complete presentation of the new and subtle geometric argument opened up by the idea of intrinsic transversality, proving linear convergence.

Intrinsic transversality, as described in Definition \ref{int}, is a geometric idea.  However, for the purposes of our proof, a variational perspective is more useful.  

Using the normalization map~
$\widehat{\mbox{~}} \colon \E\setminus \{0\} \to \E$, defined by $\widehat{z} = \frac{z}{|z|}$, we can rephrase the definition as follows.  Two sets $X,Y \subset \E$ are intrinsically transversal at a common point if there exists a constant $\kappa \in (0,1]$ (called a {\em constant of transversality}) such that nearby, any two points $x \in X \cap Y^c$ and $y \in Y \cap X^c$ have normalized displacement 
$u = \widehat{x-y}$ satisfying
\bmye \label{const}
\max \Big\{ d\big(u,N_Y(y)\big) , d\big(u,-N_X(x)\big) \Big\} ~\ge~ \kappa.
\emye

With this in mind, we have the following easy observation.

\begin{prop} \label{implication}
~~Transversality implies intrinsic transversality.  Indeed, if two closed sets $X,Y \subset \E$ intersect transversally at a point $z \in X \cap Y$, and $\theta \in (0,\pi]$ is the minimal angle between vectors in the normal cones $N_Y(z)$ and $-N_X(z)$, then intrinsic transversality holds with any constant in the interval $(0,\sin \frac{\theta}{2})$.
\end{prop}

\pf
Clearly we have
\[
\min_{|u| = 1} \max \Big\{ d\big(u,N_Y(z)\big) , d\big(u,-N_X(z)\big) \Big\} ~=~ \sin \frac{\theta}{2}.
\]
(We define $\theta = \pi$ in the trivial case when either normal cone is $\{0\}$.)
Transversality simply means $\theta > 0$, and in that case an easy limiting argument shows that inequality (\ref{const}) holds providing that $0 < \kappa < \sin \frac{\theta}{2}$.
\finpf

Several other common special cases are also covered under the umbrella of intrinsic transversality.  Examples include two closed convex sets with intersecting relative interiors, and two arbitrary intersecting sets that are polyhedra (or finite unions of polyhedra).

We mention one particular classical generalization of transversality in this context (see \cite{FW}).  Two smooth manifolds 
$X,Y \subset \E$ intersect {\em cleanly} at a point $z \in X \cap Y$ (see \cite{hormander}) if the intersection $X \cap Y$ is a manifold and the tangent spaces at $z$ satisfy
\[
T_{X \cap Y}(z) = T_{X}(z) \cap T_{Y}(z).
\]
One consequence is local linear convergence of the method of alternating projections, since, in fact, clean intersection implies intrinsic transversality.  To see this, we note that, assuming clean intersection at $z$, there exists a diffeomorphism mapping $z$ to zero and the two manifolds onto two linear subspaces \cite[Proposition C.3.1]{hormander}.  However, intrinsic transversality always holds at zero for subspaces, and an easy exercise shows that it is preserved under diffeomorphisms, so the result follows.

\section{Slope and the coupling function}\label{sec:slope}
Intrinsic transversality has its roots in a fundamental variational-analytic idea:  {\em slope}.  Denoting the extended reals $[-\infty,+\infty]$ by 
$\bR$, if a function $f \colon \E \to \bR$ is finite at a point $x \in \E$, then its slope there is
\[
|\nabla f|(x) = \limsup_{{w \to x}\atop{w \ne x}}\frac{f(x)-f(w)}{|x-w|},
\]
unless $x$ is a local minimizer, in which case we set $|\nabla f|(x) = 0$.  
The {\em limiting slope} at $x$ is
\[
\overline{|\nabla f|}(x) = \liminf_{{w \to x}\atop{f(w) \to f(x)}} |\nabla f|(w),
\]
and when $f$ is lower semicontinuous we have the following relationship:
\bmye \label{lim-eq}
|\nabla f|(x) ~\ge~ \overline{|\nabla f|}(x) ~=~ d\big(0,\partial f(x)\big),
\emye
where $\partial f$ denotes the set of subgradients \cite{VA}.  For a proof, see \cite[Propositions 1 and 2, Chapter 3]{ioffe_survey}, or, for example, \cite[Proposition~4.6]{descent_curves}.

Now consider a ``coupling'' function for the two nonempty closed sets $X,Y \subset \E$.  We define a function $\phi \colon \E^2 \to \bR$ by
\[
\phi(x,y) = \delta_X(x) + |x-y| + \delta_Y(y),
\]
where $\delta_X$ and $\delta_Y$ denote the indicator functions of the sets $X$ and and $Y$ respectively.  We also consider the marginal function $\phi_y \colon \E \to \bR$ (for any fixed point 
$y \in \E$) defined by $\phi_y(x) = \phi(x,y)$, and the marginal function 
$\phi_x \colon \E \to \bR$ (for any fixed point $x \in \E$) defined by $\phi_x(y) = \phi(x,y)$.  Standard subdifferential calculus applied to the coupling function and its marginals shows, for distinct points 
$x \in X$ and $y \in Y$, with normalized displacement $u = \widehat{x-y}$,
\[
\partial \phi_y(x) =  u + N_X(x) 
~~~\mbox{and}~~~ 
\partial \phi_x(y) = -u + N_Y(y),
\]
and furthermore
\[
\partial \phi(x,y) = \partial \phi_y(x) \times \partial \phi_x(y).
\]
Using the relationship (\ref{lim-eq}), we arrive at the following expressions for the limiting slopes:
\[
\overline{|\nabla \phi_y|}(x) = d\big(u,-N_X(x)\big)
~~~\mbox{and}~~~ 
\overline{|\nabla \phi_x|}(y) = d\big(u,N_Y(y)\big).
\]
Hence the sets $X$ and $Y$ are intrinsically transversal at a common point $z$ exactly when there exists a constant $\kappa \in (0,1]$ such that, near $z$, for any two points $x \in X \cap Y^c$ and $y \in Y \cap X^c$, the limiting slopes 
$\overline{|\nabla \phi_y|}(x)$ and  $\overline{|\nabla \phi_x|}(y)$ cannot both be strictly less than 
$\kappa$.  The following result is an interesting alternative view.

\begin{prop} \label{uniform}
Two closed sets $X,Y \subset \E$ are intrinsically transversal at a point $z \in X \cap Y$ if and only if the coupling function has slope uniformly bounded away from zero for pairs $(x,y)$ with $x \in X \cap Y^c$ and $y \in Y \cap X^c$, both near $z$.
\end{prop}

\pf
We note the relationship
\bmye \label{norm}
\big( \overline{|\nabla \phi|}(x,y) \big)^2  ~=~  \big( \overline{|\nabla \phi_y|}(x) \big)^2 + 
\big( \overline{|\nabla \phi_x|}(y) \big)^2,
\emye
from which the result follows immediately.
\finpf

\section{Error bounds and level sets}
We use the following key tool \cite[Basic Lemma, Chapter 1]{ioffe_survey}, giving a slope-based criterion for a lower semicontinuous function $f$ to have a nonempty level set
\[
[f \le \alpha] ~=~ \{ u \in \E : f(u) \le \alpha \}.
\]
We note in passing that the result (and the proof we present) holds more generally, with $\E$ replaced by an arbitrary complete metric space.

\begin{thm}[Error bound] \label{error}
~~~~Consider any lower semicontinuous function
\mbox{$f \colon \E \to \overline\R$}, finite at the point 
$x \in \E$, and constants $\delta >0$ and $\alpha < f(x)$ satisfying
\[
\inf_{w \in \E} \big\{ |\nabla f|(w) : \alpha<f(w)\le f(x),~ |w-x| \le \delta  \big\} ~>~ \frac{f(x)-\alpha}{\delta}.
\]
Then the level set $[f \le \alpha]$ is nonempty, and furthermore its distance from $x$ is no more than 
$\frac{1}{K}\big(f(x)-\alpha\big)$, where $K$ denotes the left-hand side of the inequality above. 
\end{thm}

\pf
We present two proofs.  The first, suitable in any complete metric space, applies the Ekeland variational principle to the function $g \colon \E \to \overline\R$, defined as the positive part of the function $f-\alpha$.  We deduce, for any constant 
$\gamma \in \big(\frac{1}{\delta}(f(x)-\alpha),K\big)$, the existence of a point $v$ minimizing the function
$g + \gamma|\cdot - v|$ and satisfying the properties $g(v) \le g(x)$ and 
$|v - x| \le  \frac{1}{\gamma}g(x) < \delta$.  The minimizing property of $v$ shows 
$|\nabla g|(v) \le \gamma < K$.  Hence $g(v) = 0$, since otherwise 
$|\nabla g|(v) = |\nabla f|(v) \ge K$.  The result now follows by letting $\gamma$ approach $K$.

A second, more elementary proof dispenses with the Ekeland principle when $\E$ is a Euclidean space.  Assuming as before $\frac{1}{\delta}(f(x)-\alpha) < \gamma < K$, define a set
\[
W ~=~ \{w \in \E : f(w) \le f(x) - \gamma|w-x| \big\}.
\]
The value $\beta = \inf_W f$ satisfies $\beta \le f(x)$.  Any sequence $(w_r)$ in $W$ with $f(w_r) \to \beta$ must be bounded, and any limit point $\bar w$ satisfies 
$\beta = f(\bar w) \le f(x) - \gamma|\bar w - x|$. 

We claim $\beta \le \alpha$.  If not, then $\alpha < f(\bar w) \le f(x)$ and 
\[
|\bar w - x| \le \frac{f(x) - f(\bar w)}{\gamma} < \frac{f(x) - \alpha}{\gamma} < \delta,
\]
so by definition $K \le |\nabla f|(\bar w)$.  Hence, by definition of the slope, there is a point 
$w \in \E$ satisfying
\[
f(w) < f(\bar w) - \gamma |\bar w - w| < \beta.
\]
Furthermore, the central expression is bounded above by $f(x) - \gamma(|\bar w - x|  + |\bar w - w|)$, and hence by $f(x) - \gamma |w - x|$, so $w \in W$, contradicting the definition of $\beta$.

We conclude $\beta \le \alpha$, so we have constructed a point $\bar w$ in the level set 
$[f \le \alpha]$, as required.  To see the distance estimate, notice (as in the first proof) that we can assume $f \ge \alpha$, by replacing $f$ by $\max\{f, \alpha\}$.  In that case, we have in fact proved
$\beta = \alpha$ and hence $|\bar w - x| \le \frac{1}{\gamma}(f(x) - \alpha)$.  Now let $\gamma$ approach $K$ as before.
\finpf

Our main result follows quickly from the following appealing geometric consequence of Theorem \ref{error}.  This tool concerns the distance from a point $y$ to a set $X$.  Given a trial point $x \in X$, it shows how to improve the upper bound $|x-y|$ under the assumption of a uniform lower bound on the angle between normal vectors to $X$ at points $w$ near $x$ and line segments from $w$ to $y$.

\begin{thm}[Distance decrease] \label{special}
Consider a closed set $X \subset \E$ and points $x \in X$ and $y \not\in X$ with $\rho :=|y-x|$.  Given any constant $\delta>0$, if
\[
\inf_{w \in X} \big\{ d \big( \widehat{y-w} , N_X(w) \big) : w \in B_{ \rho}(y)\cap B_{\delta}(x)  \big\}
~=~ \mu ~>~ 0,
\]
then 
\[
d(y,X) \le |y-x| - \mu\delta.
\]  
\end{thm}

\pf
We apply Theorem \ref{error} (Error bound) to the function $f = |\cdot - y| + \delta_X$.   Inequality (\ref{lim-eq}) and a repeat of our calculation of $\overline{|\nabla \phi_y|}$ in the previous section shows
\[
|\nabla f|(w) ~\ge~ \overline{|\nabla f|}(w) ~=~ d \big(\widehat{y-w},N_X(w)\big).
\]
In the notation of Theorem \ref{error}, assuming $\alpha < f(x)$, we have
\[
\mu ~\le~ K ~=~ 
\inf_{w \in X} \big\{ |\nabla f|(w) : \alpha < |w-y| \le |x-y|,~ |w-x| \le \delta \big\}.
\]
Hence, providing $\alpha$ satisfies
$\mu > \frac{1}{\delta}(|x-y| - \alpha)$, we deduce the first conclusion of Theorem \ref{error}:  the level set $[f \le \alpha]$ is nonempty, or in other words
$d(y,X) \le \alpha$.  The result now follows.
\finpf

\section{Main result}\label{sec:main_proof}
The following is our main result.
\begin{thm}[Linear convergence] \label{lin}
If two closed sets $X,Y \subset \E$ are intrinsically transversal at a point $z \in X \cap Y$, with constant $\kappa > 0$, then, for any constant $c$ in the interval $(0,\kappa)$, the method of alternating projections, initiated sufficiently near $z$, converges to a point in the intersection $X \cap Y$ with R-linear rate $1-c^2$.
\end{thm}

\pf
There exists a constant $\epsilon > 0$ such that for any distinct points 
$x \in X$ and $y \in Y$, both lying in the open ball $B_{\epsilon}(z) $, we have
\[
\max \big\{ d\big(u,N_Y(y)\big) , d\big(u,-N_X(x)\big) \big\} ~\ge~ \kappa ~~ \mbox{for} ~ u = \widehat{x-y}.
\]

Now consider any point $x \in X$, close to $z$, and $y \in P_Y(x)$.  We aim to deduce the inequality
\[
d(y,X) \le (1-c^2)|y-x|.
\]
The result will then follow by a routine induction.  We can assume $y \not\in X$, since otherwise there is nothing to prove.

Define $\rho = |x-y|$ and assume
\[
x,y \in B_{\epsilon - \kappa\rho}(z),
\]
as certainly holds if $x$ is sufficiently close to $z$.  Consider any point 
$w \in X \cap B_{\kappa\rho}(x)$.  Since $|w-x| < \kappa |x-y| \le d_Y(x)$, we see $w \not\in Y$, and in particular $w \ne y$.  We have
\[
d\big(\widehat{w-y} , N_Y(y)\big)
~\le~
d\big(\widehat{w-y} , \R_+(x-y)\big),
\]
since $x-y \in N_Y(y)$.  To bound the right-hand side, note that any nonzero vectors 
$p,q \in \E$ satisfy the inequalities
\[
d ( \widehat{p} , \R_+ q) ~\le~ 
\big| \widehat{p} - \ip{\widehat{p}}{\widehat{q}} \widehat{q} \big| ~\le~ 
\frac{|p-q|}{|q|}.
\]
Setting $p=w-y$ and $q=x-y$ shows, for any $w \in X \cap B_{\kappa\rho}(x)$, the inequality
\[
d\big(\widehat{w-y} , N_Y(y)\big) < \kappa,
\]
and hence, using intrinsic transversality
\[
d\big(\widehat{y-w} , N_X(w)\big) \ge \kappa.
\]
We next apply Theorem \ref{special} with $\delta = c\rho$.  In the notation of that result, we have 
$\mu \ge \kappa > c$, so we deduce
\[
d(y,X) \le \rho - \mu c \rho \le (1-c^2)|y-x|,
\]
as required.  The result now follows.
\finpf

Theorem \ref{basic} now follows, by applying Proposition \ref{implication}.  Indeed, if two closed sets $X,Y \subset \E$ intersect transversally at a point $z \in X \cap Y$, and the minimal angle between vectors in the closed cones $N_Y(z)$ and $-N_X(z)$ is $\theta \in (0,\pi]$,  then for any constant  
$r \in (\cos^2 \frac{\theta}{2}, 1)$, the method of alternating projections, initiated sufficiently near $z$, must converge to a point in $X \cap Y$ with R-linear rate $r$.

Another consequence of the proof technique is the following result.

\begin{thm}
Intrinsic transversality implies subtransversality.
\end{thm}

\pf
Under the assumptions of Theorem \ref{lin}, consider any point $x \in X$ near the intersection point $z$, and any point $y \in P_Y(x)$.  We proved $d(y,X) \le (1-c^2)d(x,Y)$, so by induction, alternating projections converges to a limit in the intersection $X \cap Y$ whose distance from $x$ is no more than
\[
d(x,Y) \sum_{n=0}^{\infty} (1-c^2)^n.
\]
Hence $d(x,X \cap Y) \le c^{-2} d(x,Y)$.

Now consider any point $w \in \E$ near the intersection point $z$, and any point $x \in P_X(w)$.  The triangle inequality implies
\[
d(x,Y) \le |x-w| + d(w,Y) = d(w,X) + d(w,Y).
\]
and hence
\begin{eqnarray*}
d(w,X \cap Y)	& \le & |w-x| + d(x,X \cap Y) \\
				& \le & d(w,X) + c^{-2}d(x,Y) \\
				& \le & (1+c^{-2})\big(d(w,X) + d(w,Y)\big),
\end{eqnarray*}
which proves subtransversality.
\finpf

Transversality is a strictly stronger condition than intrinsic transversality, as shown by two intersecting lines in $\R^3$.  In turn, intrinsic transversality is strictly stronger than subtransversality.  Indeed, subtransversality alone does not even guarantee local convergence of alternating projections, as the following example shows.  

Consider the two closed sets in $\R^2$
\begin{eqnarray*}
X &=& \big\{ (x,0) : x \ge 0 \big\}  \\
Y &=& \big\{(0,0)\big\} \cup \bigcup_{n=-\infty}^{\infty} 
\Big\{ \big( x,\frac{1}{2}(3^{n+1}-x)(3^{-n}x-1) \big) : 3^n \le x \le 3^{n+1} \Big\}.
\end{eqnarray*}
Geometrically, $X$ is just the nonnegative $x$-axis, and $Y$ consists of sequence of parabolic arches based on the $x$-axis and scaling linearly down towards the origin.  Subtransversality at the origin is routine to check.  However, intrinsic transversality fails:  for any integer $n$, the point 
$(2\cdot 3^n,0) \in X$ projects onto the point $(2\cdot 3^n,\frac{1}{2}3^n) \in Y$ and vice versa, so alternating projections cycles at points arbitrarily close to the origin.    

\section{Semi-algebraic intersections}\label{sec:semi}
In this last section we turn to the case where the two sets are semi-algebraic.  Alternating projections on algebraic varieties, in particular, was considered in \cite{FW}. 

Classically (and reassuringly) two smooth manifolds ``typically'' intersect trans\-versally:  almost all translations to the manifolds satisfy transversality throughout the intersection. We first show the same result for semi-algebraic sets. See \cite{tame_opt} for the role of semi-algebraic ideas in nonsmooth optimization, and for the more general ``tame'' setting, to which all results of this section extend.

\begin{thm}[Generic transversality]
Consider two closed semi-algebraic sets $X,Y \subset \E$. Then for almost every vector $e$ in  $\E$,  transversality holds at every point in the (possibly empty) intersection of the sets  $X$  and $Y-e$.
\end{thm}
\pf
Consider the set-valued mapping  $F$  from  $\E$  to $\E^2$ defined by
\[
F(z)  =  (X-z) \times (Y-z).
\]
In standard variational-analytic language \cite{ioffe_strat,ioffe_tame}, we define a value of this mapping  $(a,b) \in \E^2$ as {\em critical}  when  
$F$  is not metrically regular for  $(a,b)$  at some point 
\[
z \in F^{-1}(a,b) = (X-a) \cap (Y-b).  
\]
By a simple coderivative calculation (cf.\ \cite[p.\ 491]{alt_genproj}), this is equivalent to the sets 
$X-a$ and $Y-b$ not intersecting transversally at $z$.  The semi-algebraic Sard theorem \cite{ioffe_strat,ioffe_tame} now shows that the semi-algebraic set of critical values $(a,b) \in \E^2$ has dimension strictly less than $2\dim\E$.  The result now follows by Fubini's theorem.
\finpf

It is interesting to ask whether convergence of alternating projections is guaranteed (albeit sublinear) even without intrinsic transversality. This is particularly important, since in practice to check whether transversality holds requires knowledge of a point in the intersection. As the example in the previous section shows, without intrinsic transversality, even limit points of alternating projection iterates may not lie in the intersection. Nevertheless, such pathologies do not occur in the semi-algebraic setting. 
A key tool for establishing this result is the following theorem based on \cite[Theorem 14]{Lewis-Clarke}, a result growing out of the work of \cite{Kur} extending the classical {\L}ojasiewicz inequality for analytic functions to a nonsmooth semi-algebraic setting.  See \cite[Corollary 6.2]{tame_opt} for a related version.

\begin{thm}[Kurdyka-{\L}ojasiewicz inequality]
~~~~Given a semi-algebraic lower-semicontinuous function $f\colon\E\to\overline{\R}$, consider a bounded open set $U \subset \E$. Then, for any value $\tau\in\R$, there exists a constant $\rho >0$ and a continuous strictly increasing semi-algebraic function $\theta \colon (\tau,\tau+\rho) \to (0,+\infty)$, such that the inequality
\[
\slof (x)\geq \theta \big((f(x)\big)
\]
holds for all $x\in U$ with $\tau<f(x)<\tau+\rho$.
\end{thm}

We now arrive at the main result of this section.

\begin{thm}[Semi-algebraic intersections]\label{thm:semi-result}~  
Consider two nonempty closed semi-algebraic sets $X,Y \subset \E$, with $X$ bounded.  If the method of alternating projections starts in $Y$ and near $X$, then the distance of the iterates to the intersection $X \cap Y$ converges to zero, and hence every limit point lies in $X \cap Y$.
\end{thm}

\pf
Define a bounded open set
\[
U = \{ z : d(z,X) < 1 \}.
\]
The coupling function  $\phi$ from Section \ref{sec:slope} is semi-algebraic, so the 
Kurdyka-{\L}ojasiewicz inequality above implies 
that there exists a constant $\rho \in (0,1)$ and a continuous function  $h \colon (0,\rho) \to (0,+\infty)$ 
such that all points $x \in X$ and $y \in Y \cap U$ with $0 < |x-y| < \rho$ satisfy the inequality
\[
\overline{|\nabla\phi|}(x,y)\geq \sqrt{2}\cdot h(|x-y|).
\]
From equation~\eqref{norm}, we deduce
\begin{equation}\label{ineq:main}
\max \big\{\overline{|\nabla \phi_x|}(y), \overline{|\nabla \phi_y|}(x)\big\} ~\geq~ h(|x-y|).
\end{equation}

Suppose the initial point $y_0 \in Y$ satisfies $d(y_0,X) < \rho$.  The distance between successive iterates $|x_n - y_n|$ is decreasing:  our aim is to prove its limit is zero.  By way of contradiction, suppose in fact $|x_n-y_n| \downarrow \alpha \in (0,\rho)$.
 
Associated with the continuous function $h$,  around any nonzero vector  $v$, we  define a radially symmetric open``cusp''
\[
\mathcal{K}(v) ~:=~ \big\{z\in \E: 0<|z|< \rho ~\textrm{ and }~ 
d\big(\widehat{z},  \R_+v\big)< h\big(|z|\big)\big\},
\]
where $\R_+v$ is the ray generated by $v$. 
For each iteration $n=1,2,3,\ldots$, consider the normal vector $v_n = x_n-y_n \in N_Y(y_n)$.  
Note $|v_n| > \alpha$.  Obviously the open set 
$\ck(v_n)$ contains the point $v_n$:  we next observe that in fact it also contains a uniform neighborhood.  Specifically, we claim there exists a constant $\epsilon > 0$ such that $x-y_n \in \ck(v_n)$ whenever 
$|x-x_n| \le \epsilon$.  Otherwise there would exist a sequence $(x'_n)$ in $\E$ satisfying
\[
|x'_n - x_n| \to 0 ~~\mbox{and}~~ x'_n - y_n \not\in \ck(v_n).
\]
Notice that the distance between the vector $x'_n - y_n$ and the vector $v_n$ converges to $0$, so the same is true after we normalize them (normalization being a lipschitz map on the set 
$\{v \in \E : |v| > \frac{\alpha}{2}\}$).  But continuity of $h$ now gives the contradiction
\[
\big| \widehat{x'_n - y_n} - \widehat{v}_n \big|
~\ge~
d \big( \widehat{x'_n - y_n} , \R_+ v_n \big)
~\ge~
h\big(|x'_n - y_n|\big)
~\to~ h(\alpha)
~>~
0.
\]

We deduce that any point $x \in X$ with $|x-y_n| < \rho$ and $|x-x_n| \le \epsilon$ satisfies
\[
\overline{|\nabla \phi_{x}|}(y_n) = 
d\big(\widehat{x-y_n},N_Y(y_n)\big) \leq 
d\big(\widehat{x-y_n},\R_+v_n\}\big) < 
h\big(|x-y_n|\big),
\]
and hence, using inequality (\ref{ineq:main}), 
\begin{equation} \label{lower}
\overline{|\nabla \phi_{y_n}|}(x) \geq h\big(|x-y_{n}|\big).
\end{equation}
Denote by $\beta$ the minimum value of the strictly positive continuous function $h$ on the interval 
$\big[\alpha, d(y_0,X)\big]$, so $\beta > 0$.  Inequality (\ref{lower}) implies
\[
d \big(\widehat{y_n - x} , N_X(x) \big) \ge \beta,
\]
and hence, using Theorem~\ref{special} (Distance decrease), 
\[
|x_{n+1} - y_{n+1}| \le |x_{n+1} - y_n| \le |x_n - y_n| - \epsilon\beta,
\]
giving a contradiction for large $n$.  The result follows.
\finpf

\noindent
We remark that the technique we present here applies to much broader classes of sets --- those ``definable in an o-minimal structure'' \cite{Coste-min}.
Under a further regularity condition, \cite[Corollary 9]{noll3} proves convergence of the full alternating projection sequence in the semi-algebraic (or, more broadly, subanalytic) setting.

\bibliography{bibliography}

\def\cprime{$'$}
\begin{thebibliography}{10}

\bibitem{FW}
F.~Andersson and M.~Carlsson.
\newblock Alternating projections on nontangential manifolds.
\newblock {\em Constr. Approx.}, 38(3):489--525, 2013.

\bibitem{ABRS}
H.~Attouch, J.~Bolte, P.~Redont, and A.~Soubeyran.
\newblock Proximal alternating minimization and projection methods for
  nonconvex problems: an approach based on the {K}urdyka-\l ojasiewicz
  inequality.
\newblock {\em Math. Oper. Res.}, 35(2):438--457, 2010.

\bibitem{alt_class}
H.H. Bauschke and J.M. Borwein.
\newblock On the convergence of von {N}eumann's alternating projection
  algorithm for two sets.
\newblock {\em Set-Valued Anal.}, 1(2):185--212, 1993.

\bibitem{bapp}
H.H. Bauschke, D.R. Luke, H.M. Phan, and X.~Wang.
\newblock Restricted normal cones and the method of alternating projections:
  Applications.
\newblock {\em Set-Valued and Variational Anal.}, pages 1--27, 2013.

\bibitem{btheory}
H.H. Bauschke, D.R. Luke, H.M. Phan, and X.~Wang.
\newblock Restricted normal cones and the method of alternating projections:
  Theory.
\newblock {\em Set-Valued and Variational Anal.}, pages 1--43, 2013.

\bibitem{BCR}
J.~Bochnak, M.~Coste, and M.-F. Roy.
\newblock {\em Real algebraic geometry}, volume~36 of {\em Ergebnisse der
  Mathematik und ihrer Grenzgebiete (3) [Results in Mathematics and Related
  Areas (3)]}.
\newblock Springer-Verlag, Berlin, 1998.
\newblock Translated from the 1987 French original, Revised by the authors.

\bibitem{Lewis-Clarke}
J.~Bolte, A.~Daniilidis, A.S. Lewis, and M.~Shiota.
\newblock {C}larke subgradients of stratifiable functions.
\newblock {\em SIAM J. Optim.}, 18(2):556--572, 2007.

\bibitem{Borwein-Zhu}
J.M. Borwein and Q.J. Zhu.
\newblock {\em Techniques of Variational Analysis}.
\newblock Springer Verlag, New York, 2005.

\bibitem{conv_alt}
L.M. Bregman.
\newblock The method of successive projection for finding a common point of
  convex sets.
\newblock {\em Sov. Math. Dokl}, 6:688--692, 1965.

\bibitem{CLSW}
F.H. Clarke, Yu. Ledyaev, R.I. Stern, and P.R. Wolenski.
\newblock {\em Nonsmooth Analysis and Control Theory}.
\newblock Texts in Math. 178, Springer, New York, 1998.

\bibitem{Coste-min}
M.~Coste.
\newblock {\em An introduction to o-minimal geometry}.
\newblock RAAG Notes, 81 pages, Institut de Recherche Math\'{e}matiques de
  Rennes, November 1999.

\bibitem{Coste-semi}
M.~Coste.
\newblock {\em An {I}ntroduction to {S}emialgebraic {G}eometry}.
\newblock RAAG Notes, 78 pages, Institut de Recherche Math\'{e}matiques de
  Rennes, October 2002.

\bibitem{dima-thesis}
D.~Drusvyatskiy.
\newblock {\em Slope and Geometry in Variational Mathematics}.
\newblock PhD thesis, Cornell University, 2013.

\bibitem{descent_curves}
D.~Drusvyatskiy, A.D. Ioffe, and A.S. Lewis.
\newblock Curves of descent.
\newblock {\em SIAM J. Control Optim.}, 53(1):114--138, 2015.

\bibitem{coupling}
D.~Drusvyatskiy, A.D. Ioffe, and A.S. Lewis.
\newblock Alternating projections and coupling slope.
\newblock {\em arxiv:1401.7569v1}, January 2014.

\bibitem{GPR}
L.G. Gubin, B.T. Polyak, and E.V. Raik.
\newblock The method of projections for finding the common point of convex
  sets.
\newblock {\em USSR Computational Mathematics and Mathematical Physics},
  7(6):1--24, 1967.

\bibitem{alt_glob}
R.~Hesse and D.R. Luke.
\newblock Nonconvex notions of regularity and convergence of fundamental
  algorithms for feasibility problems.
\newblock {\em Under review, arXiv:1212.3349 [math.OC]}, 2012.

\bibitem{hormander}
L.~H{\"o}rmander.
\newblock {\em The analysis of linear partial differential operators. {III}}.
\newblock Classics in Mathematics. Springer, Berlin, 2007.
\newblock Pseudo-differential operators, Reprint of the 1994 edition.

\bibitem{ioffe_survey}
A.D. Ioffe.
\newblock Metric regularity and subdifferential calculus.
\newblock {\em Uspekhi Mat. Nauk}, 55(3(333)):103--162, 2000.

\bibitem{ioffe_tame}
A.D. Ioffe.
\newblock A {S}ard theorem for tame set-valued mappings.
\newblock {\em J. Math. Anal. Appl.}, 335(2):882--901, 2007.

\bibitem{ioffe_strat}
A.D. Ioffe.
\newblock Critical values of set-valued maps with stratifiable graphs.
  {E}xtensions of {S}ard and {S}male-{S}ard theorems.
\newblock {\em Proc. Amer. Math. Soc.}, 136(9):3111--3119, 2008.

\bibitem{tame_opt}
A.D. Ioffe.
\newblock An invitation to tame optimization.
\newblock {\em SIAM J. Optim.}, 19(4):1894--1917, 2009.

\bibitem{ioffe-survey}
A.D. Ioffe.
\newblock Metric regularity. theory and applications -- a survey.
\newblock {\em Preprint arxiv:1505.07920}, 2015.

\bibitem{kruger}
A.Y. Kruger and N.H. Thao.
\newblock Quantitative characterizations of regularity properties of
  collections of sets.
\newblock {\em J. Optim. Theory Appl.}, 164(1):41--67, 2015.

\bibitem{kruger-thao}
A.Y. Kruger and N.H. Thao.
\newblock Regularity of collections of sets and convergence of inexact
  alternating projections.
\newblock {\em arXiv:1501.04191}, 2015.

\bibitem{Kur}
K.~Kurdyka.
\newblock On gradients of functions definable in o-minimal structures.
\newblock {\em Ann. Inst. Fourier (Grenoble)}, 48(3):769--783, 1998.

\bibitem{L-ICM}
A.S. Lewis.
\newblock {\em Nonsmooth optimization: conditioning, convergence and
  semialgebraic models}.
\newblock In S.Y. Jang, Y.R. Kim, D.-W. Lee, and I. Yie, editors, Proceedings
  of the International Congress of Mathematicians, Seoul, volume IV: Invited
  Lectures, pages 872-895, Seoul, Korea, 2014.

\bibitem{alt_genproj}
A.S. Lewis, D.R. Luke, and J.~Malick.
\newblock Local linear convergence for alternating and averaged nonconvex
  projections.
\newblock {\em Found. Comput. Math.}, 9(4):485--513, 2009.

\bibitem{alt_man}
A.S. Lewis and J.~Malick.
\newblock Alternating projections on manifolds.
\newblock {\em Math. Oper. Res.}, 33(1):216--234, 2008.

\bibitem{Mord_1}
B.S. Mordukhovich.
\newblock {\em Variational Analysis and Generalized Differentiation I: Basic
  Theory}.
\newblock Grundlehren der mathematischen Wissenschaften, Vol 330, Springer,
  Berlin, 2006.

\bibitem{Mord_2}
B.S. Mordukhovich.
\newblock {\em Variational Analysis and Generalized Differentiation II:
  Applications}.
\newblock Grundlehren der mathematischen Wissenschaften, Vol 331, Springer,
  Berlin, 2006.

\bibitem{noll3}
D.~Noll and A.~Rondepierre.
\newblock On local convergence of the method of alternating projections.
\newblock {\em Found. Comput. Math.}, 16(2):425--455, 2016.

\bibitem{noll}
D.~Noll and A.~Rondepierre.
\newblock On local convergence of the method of alternating projections.
\newblock {\em arXiv:1312.5681v1 [math.OC]}, December 2013.

\bibitem{noll2}
D.~Noll and A.~Rondepierre.
\newblock On local convergence of the method of alternating projections.
\newblock {\em arXiv:1312.5681v2 [math.OC]}, September 2014.

\bibitem{VA}
R.T. Rockafellar and R.J-B. Wets.
\newblock {\em Variational {A}nalysis}.
\newblock Grundlehren der mathematischen Wissenschaften, Vol 317, Springer,
  Berlin, 1998.

\bibitem{N_proj}
J.~von Neumann.
\newblock {\em Functional {O}perators. {II}. {T}he {G}eometry of {O}rthogonal
  {S}paces}.
\newblock Annals of Mathematics Studies, no. 22. Princeton University Press,
  Princeton, N. J., 1950.

\end{thebibliography}
\bibliographystyle{plain}

\end{document}